\documentclass[a4paper,reqno]{amsart}
\usepackage{amssymb, amsmath, amscd}
\usepackage{color}

\newtheorem{theorem}{Theorem}[section]

\newtheorem{remark}[theorem]{Remark}

\def\qedbox{\hbox{$\rlap{$\sqcap$}\sqcup$}}
\def\WW{\mathcal{W}}\def\WWW{\tilde{\mathcal{W}}}

\makeatletter
 \@addtoreset{equation}{section}
\makeatother
\begin{document}
\title[Geometric realizations]
 {Geometric realizations of Kaehler and of para-Kaehler curvature models}
\author[M. Brozos-V\'azquez et. al.]{M. Brozos-V\'azquez,
P. Gilkey, and E. Merino}
\address{MBV: E. U. P. Ferrol, University of A Coru\~na, Spain\\
E-mail: mbrozos@udc.es}
\address{PG: Mathematics Department, University of Oregon\\
  Eugene OR 97403 USA\\
  E-mail: gilkey@uoregon.edu}
\address{EM: E. P. S. Ferrol, University of A Coru\~na, Spain\\
E-mail: uxiomer@cdf.udc.es}

\begin{abstract}{We show that every Kaehler algebraic curvature tensor
is geometrically realizable by a Kaehler manifold of constant scalar curvature.  We also show
that every para-Kaehler algebraic curvature tensor is geometrically realizable by a para-Kaehler
manifold of constant scalar curvature.}\end{abstract}
\maketitle
\section{Introduction}

Curvature is a central aspect of modern differential geometry. One can relate curvature to the
underlying geometry and topology of a manifold
\cite{BGNS08,DJ08,T08}, examine analytic properties which are influenced by curvature \cite{K08}, and study purely
algebraic properties of curvature
\cite{GHVV08,GN07,MP09}.
The study of Hermitian and Kaehler geometry is an active current field of investigation
that continues to be important both in pure mathematics \cite{B08} and in mathematical physics \cite{B07,CM08,W08}.
Furthermore, the study of para-Kaehler geometry, which is a neutral signature analogue,
also is an active research area
\cite{CVLS06,CDG08,DDGMMV07}.

In this paper, we extend previous investigations \cite{BGKSW09,BGSV09,GN08x,GNW09x,GNW09} to discuss questions of geometric
 realizability -- when can an algebraic curvature tensor, which is a purely algebraic object, be  realized
by a Riemannian manifold in suitable contexts and what are the resulting geometric  constraints, if any. We will
focus our attention on the Kaehler and para-Kaehler settings and give necessary and sufficient linear conditions on the
curvature to ensure that a given curvature model is geometrically  realizable by a Kaehler manifold or by a
para-Kaehler manifold. Imposing the condition that the manifold in question has constant scalar curvature yields no
additional restrictions.

\subsection{Hermitian manifolds}
We begin with a brief review of previously known results. Let $\mathcal{M}:=(M,g)$ be a Riemannian manifold of dimension
$m$; we shall always assume that
$m=2\bar m$ is even and that $m\ge4$. Let $R$ be the curvature tensor of
$\mathcal{M}$; $R$ satisfies the following identity for all tangent vectors
 $x,y,z,w$:
\begin{equation}\label{eqn-1.a}
\begin{array}{l}
R(x,y,z,w)=-R(y,x,z,w)=R(z,w,x,y),\\
0=R(x,y,z,w)+R(y,z,x,w)+R(z,x,y,w)\,.\vphantom{\vrule height 11pt}
\end{array}\end{equation}
Suppose there exists an {\it almost complex structure} $J$ on $M$, i.e. an endomorphism
$J$ of the tangent bundle
$TM$ so that
$J^2=-\operatorname{id}$. We also assume that $J$ is {\it Hermitian}, i.e. that $J^*g=g$; in this setting, the triple
$\mathcal{C}:=(M,g,J)$ is said to be an {\it almost Hermitian manifold}.
An almost Hermitian manifold $\mathcal{C}$ is said to be a {\it
Hermitian manifold} if
$J$ is {\it integrable}; this means that  the {\it Nijenhuis tensor}
$$N_{J}(x,y):=[x,y]+J[Jx,y]
 +J[x,Jy]-[Jx,Jy]$$
vanishes
or, equivalently, by the Newlander--Nirenberg Theorem \cite{NN57} that every point of $M$ has a neighborhood
with local coordinates
 $(x_1,y_1,...,x_{\bar m},y_{\bar m})$ so that
$$J\partial_{x_i}=\partial_{y_i}\quad\text{and}
\quad J\partial_{y_i}=-\partial_{x_i}\,.$$

If $\mathcal{C}$ is a Hermitian manifold, then the Riemann curvature tensor satisfies an extra identity
discovered by Gray
\cite{gray}:
\begin{theorem}\label{thm-1.1}
If $\mathcal{C}$ is a Hermitian manifold, then
\begin{eqnarray}
0&=&R(x,y,z,w)+R(Jx,Jy,Jz,Jw)-R(Jx,Jy,z,w)\nonumber\\
&-&R(Jx,y,J z,w)-R(Jx,y,z,Jw)-R(x,Jy,Jz,w)\label{eqn-1.b}\\
&-&R(x,Jy,z,Jw)-R(x,y,Jz,Jw)\quad\forall\quad x,\ y,\ z,\ w\nonumber\,.
\end{eqnarray}
\end{theorem}

\begin{remark}\label{rmk-1.2}
\rm Theorem \ref{thm-1.1} shows that the integrability of the almost complex structure implies a relation in the
curvature. Let $\{x_1,y_1,x_2,y_2,x_3,y_3\}$ be coordinates on $\mathbb{R}^6$ with the standard flat metric
$ds^2=dx_1^2+dy_1^2+dx_2^2+dy_2^2+dx_3^2+dy_3^2$. Define
$$\begin{array}{ll}
J\partial_{x_1}=\phantom{-}\cos(x_3)\partial_{y_1}+\sin(x_3)\partial_{y_2},&
J\partial_{y_1}=-\cos(x_3)\partial_{x_1}+\sin(x_3)\partial_{x_2},\\
J\partial_{x_2}=-\sin(x_3)\partial_{y_1}+\cos(x_3)\partial_{y_2},&
J\partial_{y_2}=-\sin(x_3)\partial_{x_1}-\cos(x_3)\partial_{x_2},\\
J\partial_{x_3}=\partial_{y_3},&J\partial_{y_3}=-\partial_{x_3}\,.
\end{array}$$
The resulting structure is a flat almost Hermitian manifold which is not integrable. Thus relations in curvature do not imply
integrability.
\end{remark}

\subsection{Hermitian curvature models} We now pass to the algebraic setting. Let $V$ be a real vector space of
dimension
$m$ which is equipped with a positive definite inner product $\langle\cdot,\cdot\rangle$. We say that
$A\in\otimes^4V^*$ is an {\it algebraic curvature tensor} if $A$ has the symmetries of Equation (\ref{eqn-1.a});
let
$\mathfrak{A}=\mathfrak{A}(V)$ be the subspace of $\otimes^4V^*$ which consists of all such tensors. Fix
$A\in\mathfrak{A}$ and let $\mathfrak{M}:=(V,\langle\cdot,\cdot\rangle,A)$ be the associated {\it curvature model}.
We
say
$\mathfrak{C}:=(V,\langle\cdot,\cdot\rangle,J,A)$ is an {\it almost Hermitian curvature model} if
$J$ is a linear map of $V$ with $J^2=-\operatorname{id}$ and
$J^*\langle\cdot,\cdot\rangle=\langle\cdot,\cdot\rangle$.
Let $\{e_i\}$ be an orthonormal basis for $(V,\langle\cdot,\cdot\rangle)$. We adopt the {\it Einstein} convention
and sum over repeated indices and define the {\it Ricci tensor}
$\rho$, the {\it
$\star$-Ricci tensor}
$\rho^\star$, the {\it scalar curvature} $\tau$, and the {\it $\star$-scalar  curvature} $\tau^\star$ by contracting
indices:
$$\begin{array}{ll}
\rho(x,y):=A(x,e_i,e_i,y),&
 \tau:=A(e_i,e_j,e_j,e_i),\\
 \rho^\star(x,y):= A(x,e_i,Je_i,Jy),
 &\tau^\star:=A(e_i,e_j,Je_j,Je_i)\,.\vphantom{\vrule height 12pt}
\end{array}$$

If in addition the Gray identity given in Equation (\ref{eqn-1.b}) is satisfied, then
$\mathfrak{C}$ is said to be a {\it Hermitian curvature model}. We say that
$\mathfrak{C}$ is {\it geometrically realizable} by an almost Hermitian manifold $\mathcal{C}=(M,g,J)$ if
there exists a point $P$ of $M$ and an isomorphism $\phi:T_PM\rightarrow V$ so that
$\phi^*\langle\cdot,\cdot\rangle=g_P$, $\phi^*J=J_P$, and $\phi^*A=R_P$. We refer to
 \cite{BGKSW09,BGKS09} for the proof of the following result which provides a converse to Theorem \ref{thm-1.1}:

\begin{theorem}\label{thm-1.3}
\
\begin{enumerate}
\item
If $\mathfrak{C}$ is an almost Hermitian curvature model,  then $\mathfrak{C}$ is geometrically realizable by an almostHermitian manifold with constant scalar curvature  and with constant $\star$-scalar curvature.
\item If $\mathfrak{C}$ is a Hermitian curvature model, then $\mathfrak{C}$ is geometrically realizable by a Hermitian
manifold with constant scalar curvature and with constant $\star$-scalar curvature.
\end{enumerate}
\end{theorem}

\begin{remark}\label{thm-1.4}
\rm Theorem \ref{thm-1.3} shows that the existence of an almost Hermitian structure imposes no
additional relations on the curvature other than those generated by Equation (\ref{eqn-1.a}).
Furthermore, Equations (\ref{eqn-1.a}) and (\ref{eqn-1.b}) generate the universal symmetries of the
curvature tensor of a Hermitian manifold. Finally, assuming that the scalar curvature and $\star$-scalar
curvature are constant  imposes no additional relations.
\end{remark}

\subsection{Kaehler geometry}The {\it Kaehler} form of an almost Hermitian manifold
$\mathcal{C}$ is defined by setting $\Omega(x,y):=g(x,Jy)$. One says that $\mathcal{C}$ is a {\it Kaehler} manifold if
 $J$ is integrable and $d\Omega=0$ or, equivalently, if $\nabla J=0$. One has:
\begin{theorem}\label{thm-1.5}
If $\mathcal{C}$ is a Kaehler manifold, then:
\begin{equation}\label{eqn-1.c}
R(x,y,z,w)=R(Jx,Jy,z,w)\quad\forall\quad x,\ y,\ z, \ w\,.
\end{equation}\end{theorem}
An almost Hermitian curvature model
$\mathfrak{C}=(V,\langle\cdot,\cdot\rangle,J,A)$ is said to be a {\it Kaehler curvature model} if
Equation (\ref{eqn-1.c}) is satisfied; this necessarily implies Equation (\ref{eqn-1.b}) is satisfied
so any Kaehler curvature model is a Hermitian curvature model. In Section
\ref{sect-2}, we shall establish the following result which shows that Equations (\ref{eqn-1.a}) and
(\ref{eqn-1.c}) generate the universal symmetries of the curvature tensor of a Kaehler manifold and which
is a converse to Theorem \ref{thm-1.5}. In this setting, necessarily $\tau=\tau^\star$. In Section \ref{sect-4}, we
show that the realization can also be taken to have $\tau$ constant.
\begin{theorem}\label{thm-1.6}
Let $\mathfrak{C}$ be a Kaehler curvature model.  Then $\mathfrak{C}$ is geometrically
realizable by a Kaehler manifold of constant scalar curvature.
\end{theorem}

\begin{remark}\label{rmk-1.7}
\rm The methods used in  \cite{BGSV09} can be used to extend Theorem
\ref{thm-1.6} to the indefinite setting; we omit details in the interests of brevity.
\end{remark}

\begin{remark}\label{rmk-1.8}
\rm Theorems \ref{thm-1.1} and \ref{thm-1.3} show that Equation (\ref{eqn-1.b}) provides necessary and sufficient linear identities for a Hermitian curvature model to be geometrically realizable by a Hermitian manifolds. Similarly, Theorems \ref{thm-1.5} and \ref{thm-1.6} show that Equation (\ref{eqn-1.c}) provides necessary and sufficient
linear identities for a curvature model to be geometrically realizable by a
Kaehler manifold. There are examples where one has relations rather than
identities. For example,  one says that an almost Hermitian manifold is {\it
almost Kaehler} if
$d\Omega=0$. In this setting, we have $\tau^\star-\tau=\frac12\vert\nabla J\vert^2$ and thus the curvature lies in the
half-space defined by the relation
$\tau^\star\geq\tau$. This shows that an almost Hermitian curvature model with $\tau>\tau^\star$ is not geometrically realizable by an almost Kaehler manifold. We refer to
\cite{AD01,DDGMMV07} for further details concerning almost Kaehler
manifolds in both the Riemannian and the higher signature settings.
\end{remark}

\subsection{Para-Hermitian manifolds}
We shall say that the triple
$\tilde{\mathcal{C}}:=(M,g,\tilde J)$ is an {\it almost para-Hermitian manifold} if $g$ is a pseudo-Riemannian
metric on $M$ of neutral signature $(\bar m,\bar m)$ and if
$\tilde J$ is a linear map of $TM$ satisfying
$$\tilde J^2=\operatorname{id}\quad\text{and}\quad
  \tilde J^*g=-g\,.$$

If $\tilde{\mathcal{C}}$ is an almost para-Hermitian manifold, then one says that
$\tilde J$ is {\it integrable} if the {\it para-Nijenhuis tensor}
$$N_{\tilde{J}}(x,y):=[x,y]-\tilde{J}[\tilde{J}x,y]
-\tilde{J}[x,\tilde{J}y]+[\tilde{J}x,\tilde{J}y]$$
vanishes or, equivalently, every point of $M$ has a neighborhood with local coordinates  $(x_1,y_1,...,x_{\bar m},y_{\bar m})$
so that
$$\tilde{J}\partial_{x_i}=\partial_{y_i}\quad\text{and}\quad\tilde{J}\partial_{y_i}=\partial_{x_i}\,.$$
 If $\tilde{J}$ is integrable then $\tilde{\mathcal{C}}$ is called a {\it para-Hermitian manifold}. Theorem \ref{thm-1.1} generalizes to this setting \cite{BGSV09} to become the following result -- note the changes in
sign from Equation (\ref{eqn-1.b}):
\begin{theorem}\label{thm-1.9}
If $\tilde{\mathcal{C}}$ is a para-Hermitian
manifold, then
\begin{eqnarray}
0&=&R(x,y,z,w)+R(\tilde Jx,\tilde Jy,\tilde Jz,\tilde Jw)+R(\tilde Jx,\tilde
Jy,z,w)\nonumber\\
&+&R(\tilde Jx,y,\tilde Jz,w)+R(\tilde Jx,y,z,\tilde Jw)+R(x,\tilde Jy,\tilde Jz,w)\label{eqn-1.d}\\
&+&R(x,\tilde Jy,z,\tilde Jw)+R(x,y,\tilde Jz,\tilde
Jw)\quad\forall\quad x,\ y,\ z,\ w\,.\nonumber
\end{eqnarray}\end{theorem}

\subsection{Para-Hermitian curvature models} One defines the notion of an almost para-Hermitian curvature model
$\tilde{\mathfrak{C}}:=(V,\langle\cdot,\cdot\rangle,\tilde J,A)$ similarly; if $\tilde{\mathfrak{C}}$
satisfies the relations of Equation (\ref{eqn-1.d}), then
$\tilde{\mathfrak{C}}$ is said to be a para-Hermitian curvature model. Theorem \ref{thm-1.3} extends to this setting
\cite{BGKSW09,BGSV09}:
\begin{theorem}\label{thm-1.10}
\
\begin{enumerate}
\item Let $\tilde{\mathfrak{C}}$ be an almost para-Hermitian curvature model. Then $\tilde{\mathfrak{C}}$ is
geometrically realizable by an almost para-Hermitian manifold with constant scalar curvature and with constant
$\star$-scalar curvature.
\item Let $\tilde{\mathfrak{C}}$ be a para-Hermitian curvature model. Then $\mathfrak{C}$ is geometrically realizable
by a para-Hermitian manifold with constant scalar curvature and with constant $\star$-scalar curvature.
\end{enumerate}\end{theorem}

\subsection{Para-Kaehler geometry} One defines the {\it para-Kaehler}  form of
an almost para-Hermitian manifold
$\tilde{\mathcal{C}}=(M,g,\tilde J)$ by setting
$\tilde\Omega(x,y):=g(x,Jy)$. We say that
$\tilde{\mathcal{C}}$ is {\it para-Kaehler} if  $\tilde J$ is integrable and $d\tilde\Omega=0$ or, equivalently, if
$\nabla\tilde J=0$. In this setting one has (note the change in sign from Equation (\ref{eqn-1.c})):
\begin{theorem}\label{thm-1.11}
If $\tilde{\mathcal{C}}$ is a para-Kaehler manifold, then:
\begin{equation}\label{eqn-1.e}
R(x,y,z,w)=-R(\tilde Jx,\tilde Jy,z,w)\quad\forall\quad x,\ y,\ z,\ w\,.
\end{equation}\end{theorem}
We say an almost para-Hermitian curvature model $\tilde{\mathfrak{C}}$ is a {\it para-Kaehler} curvature model if the
relations of Equation (\ref{eqn-1.e}) hold; this implies the relations of Equation (\ref{eqn-1.d})
hold and thus $\tilde{\mathfrak{C}}$ is also a para-Hermitian curvature model. Theorem \ref{thm-1.6} generalizes to this
setting to become:

\begin{theorem}\label{thm-1.12}
If $\tilde{\mathfrak{C}}$ is a para-Kaehler curvature model, then
$\tilde{\mathfrak{C}}$ is geometrically realizable by a para-Kaehler manifold of constant scalar curvature.
\end{theorem}

\subsection{Outline of the paper} We show in Section \ref{sect-2} (resp. in Section
\ref{sect-3}) that any Kaehler (resp. para-Kaehler) curvature model can be geometrically realized by a
Kaehler (resp. para-Kaehler) manifold. In Section \ref{sect-4}, we show the realizations can be chosen to have
constant scalar curvature.

The decomposition of the space of algebraic curvature tensors under the action of the unitary group was given by
Tricerri and Vanhecke \cite{TV81} and is summarized in Theorem \ref{thm-2.1}; it plays a central role in the analysis of
Section
\ref{sect-2} -- a similar analysis in the para-Kaehler setting is performed in Section \ref{sect-3}.
The Cauchy-Kovalevskaya Theorem
is used in Section \ref{sect-4} to show that the Kaehler and para-Kaehler realizations in question can be chosen to have
constant scalar curvature.
\section{Kaehler Geometry}\label{sect-2}

\subsection{Curvature decomposition under the unitary group} Let $(V,\langle\cdot,\cdot\rangle,J)$ be a Hermitian
structure. Let
$\{e_i\}$ be an orthonormal basis for
$(V,\langle\cdot,\cdot\rangle)$. Let $\mathcal{U}$ be the unitary
group:
$$
\mathcal{U}:=
\{T\in\operatorname{GL}(V):TJ=JT\quad\text{and}\quad
T^*\langle\cdot,\cdot\rangle=\langle\cdot,\cdot\rangle\}\,.
$$
Define elements $\pi_1$ and $\pi_2$ of $\mathfrak{A}$ by setting:
\begin{eqnarray*}
&&\pi_1(x,y,z,w):=\langle x,w\rangle\langle y,z\rangle-\langle x,z\rangle\langle y,w\rangle,\\&&\pi_2(x,y,z,w):=\langle Jx,w\rangle\langle Jy,z\rangle-\langle Jx,z\rangle\langle Jy,w\rangle
  -2\langle Jx,y\rangle\langle Jz,w\rangle\,.
\end{eqnarray*}
 Let $S^2(V^*)$ denote the set of symmetric $2$-tensors and let $\theta\in S^2(V^*)$ with $J^*\theta=\theta$. Define elements  $\phi (\theta)$ and $\psi(\theta)$ of $\mathfrak{A}(V)$
by setting:
$$\begin{array}{l}
(\phi(\theta))(x,y,z,w):=\langle x,w\rangle\theta(y,z)-\langle x,z\rangle\theta(y,w)\\
\qquad\qquad\qquad\quad+\theta(x,w)\langle y,z\rangle-\theta(x,z)\langle y,w\rangle,
\vphantom{\vrule height 10pt}\\
(\psi(\theta))(x,y,z,w):=\langle Jx,w\rangle\theta(Jy,z)-\langle Jx,z\rangle\theta(Jy,w)
  -2\langle Jx,y\rangle\theta(Jz,w)\vphantom{\vrule height 10pt}\\
\qquad\qquad\qquad\quad+\theta(Jx,w)\langle Jy,z\rangle-\theta(Jx,z)\langle Jy,w\rangle
  -2\theta(Jx,y)\langle Jz,w\rangle\,.\vphantom{\vrule height 10pt}
\end{array}$$
The following result is
due to Tricerri and Vanhecke
\cite{TV81} in the Riemannian setting -- the extension to the higher signature context is not difficult
\cite{BGSV09}:
\begin{theorem}\label{thm-2.1}
Let $(V,\langle\cdot,\cdot\rangle,J)$ be a Hermitian structure.
\begin{enumerate}
\item We have the following orthogonal direct sum decomposition of $\mathfrak{A}$ into irreducible $\mathcal{U}$ modules:
\begin{enumerate}
\item If $m=4$,
$\mathfrak{A}=\mathcal{W}_1\oplus\mathcal{W}_2\oplus\mathcal{W}_3\oplus\mathcal{W}_4\oplus\mathcal{W}_7
\oplus\mathcal{W}_8\oplus\mathcal{W}_9$.
\item If $m=6$,
$\mathfrak{A}=\mathcal{W}_1\oplus\mathcal{W}_2\oplus\mathcal{W}_3\oplus\mathcal{W}_4\oplus\mathcal{W}_5\oplus\mathcal{W}_7
\oplus\mathcal{W}_8\oplus\mathcal{W}_9\oplus\mathcal{W}_{10}$.
\item If $m\ge8$,
$\mathfrak{A}=\mathcal{W}_1\oplus\mathcal{W}_2\oplus\mathcal{W}_3\oplus\mathcal{W}_4\oplus\mathcal{W}_5\oplus\mathcal{W}_6\oplus\mathcal{W}_7
\oplus\mathcal{W}_8\oplus\mathcal{W}_9\oplus\mathcal{W}_{10}$.
\end{enumerate}
We have $\mathcal{W}_1\approx\mathcal{W}_4$ and, if $m\ge6$, $\mathcal{W}_2\approx\mathcal{W}_5$. The other $\mathcal{U}$ modules
appear with multiplicity 1.
\item  $\mathfrak{C}=(V,\langle\cdot,\cdot\rangle,J,A)$ is a Kaehler model if
and only if
$A\in\WW_1\oplus\WW_2\oplus\WW_3$.
\item If $A\in\WW_1\oplus\WW_2\oplus\WW_3$, then the projections $p_i$ of $A$ on $\WW_i$ are given by:
\begin{enumerate}
\item $p_1A= \frac1{4\bar m(\bar m+1)}(\pi_1+\pi_2)\tau$.
\item $p_2A= \frac1{4(\bar m+2)}(\phi+\psi)(2\rho-\frac1{\bar m}\tau\langle\cdot,\cdot\rangle)$.
\item $p_3A=A-p_1A-p_2A$.
\end{enumerate}\end{enumerate}
\end{theorem}

\subsection{Realizability of Kaehler curvature models}
Let
$\{u_1,...,u_m\}$ be the canonical coordinates on $\mathbb{R}^m$.  Set $\partial_i:=\frac{\partial}{\partial
u_i}$ .
Let $J$ be the  canonical integrable almost complex structure on $\mathbb{R}^m$ given by:
\begin{equation}\label{eqn-2.a}
J\partial_{1+2k}=\partial_{2+2k}\quad\text{and}\quad
J\partial_{2+2k}=-\partial_{1+2k}\text{ for }0\le k<\bar m\,.
\end{equation}
Let $g_0$ be the usual flat Hermitian metric on $\mathbb{R}^m$:
$$g_{0,ij}:=\left\{\begin{array}{rl}
1&\text{if }i=j,\\
0&\text{otherwise}\,.\end{array}\right.$$
We identify $(V,\langle\cdot,\cdot\rangle,J)$ with $(\mathbb{R}^m,g_0,J)$. Let $S_+^2(V^*)$ be the set of all
symmetric $2$-tensors $\theta$ such that $J^*\theta=\theta$.  If $\Theta\in S^2_+(V^*)\otimes S^2(V^*)$, then set:
\begin{equation}\label{eqn-2.b}
g_{\Theta,ij}:=g_{0,ij}+\Theta_{ijkl}u^ku^l\,.
\end{equation}
This is positive definite at the origin and hence defines a Riemannian metric which is invariant under $J$
on some neighborhood of the origin. Let $\mathcal{L}\Theta\in\mathfrak{A}(V)$
be the curvature of the metric $g_\Theta$ at the origin;
$$(\mathcal{L}\Theta)(x,y,z,w)= \Theta(x,z,y,w)+\Theta(y,w,x,z)-\Theta(x,w,y,z)-\Theta(y,z,x,w)\,.$$

The map
$\Theta\rightarrow d\Omega_{g_\Theta}$ defines a linear map
$$K_J:S^2_+(V^*)\otimes S^2(V^*)\rightarrow\Lambda^3(V^*)\otimes V^*$$
which is given by
\begin{equation}\label{eqn-2.c}
\{(K_J\Theta)(x,y,z)\}(w):=\Theta(x,Jy,z,w)
+\Theta(y,Jz,x,w)+\Theta(z,Jx,y,w)\,.
\end{equation}
 This shows that $\ker(K_J)$ is invariant under the action of $\mathcal{U}$.
Clearly $\Theta\in\ker(K_J)$ if and only if $g_\Theta$ is a Kaehler metric;
$$\mathcal{L}:\ker(K_J)\rightarrow W_1\oplus W_2\oplus W_3$$
is a linear map which is equivariant with respect to the action $\mathcal{U}$.
To show every Kaehler curvature model is geometrically realizable by a Kaehler metric, it suffices to show that
$\mathcal{L}$ is surjective. Take
$$\Theta=\frac12(e^1\otimes e^1+e^2\otimes e^2)\otimes (e^1\otimes e^1+e^2\otimes e^2)$$
so that the metric has the form
$$g_\Theta=(du_1^2+...+du_{m}^2)+\frac12
(u_1^2+u_2^2)(du_1^2+du_2^2)\,.$$
The metric $g_\Theta$ is Kaehler since it takes the form
$M_2\times\mathbb{C}^{\bar m-1}$ where $M_2$ is a Riemann surface. Thus $\Theta\in\ker(K_J)$. Furthermore, the only
non-zero curvature components of the curvature tensor at the origin, up to the usual $\mathbb{Z}_2$ symmetries, are
$$R(e_1,e_2,e_2,e_1)=1\,.$$

The non-zero components of $\rho$ are $\rho(e_1,e_1)= \rho(e_2,e_2)=1$. We compute:
$$\begin{array}{ll}
 \pi_1(e_3,e_4,e_4,e_3)=1,& \pi_2(e_3,e_4,e_4,e_3)=3,\\
\tau=2,& p_1A(e_3,e_4,e_4,e_3)=\frac2{\bar m(\bar m+1)}\,.
\end{array}$$
Thus the component of $A$ in $\WW_1$ is non-zero. Similarly, we compute:
$$\begin{array}{ll}
(\phi(\rho))(e_3,e_4,e_4,e_3)=0,& (\psi(\rho))(e_3,e_4,e_4,e_3)=0,\\
(\phi(\langle\cdot,\cdot\rangle))(e_3,e_4,e_4,e_3)=2,&(\psi(\langle\cdot,\cdot\rangle))(e_3,e_4,e_4,e_3)=6,\\
 (p_2A)(e_3,e_4,e_4,e_3)= -\frac4{\bar m(\bar m+2)}\,.
\end{array}$$
This shows the component of $A$ in $\WW_2$ is non-zero. We have
$$((\operatorname{id}-p_1 -p_2)A)(e_3,e_4,e_4,e_3)=-\textstyle\frac{2(\bar m+2)-4(\bar m+1)}{\bar m(\bar m+1)(\bar m+2)}=-\frac2{(\bar m+1)(\bar m+2)}$$
and thus the component of $A$ in $W_3$ is non-zero. Thus $A$ has non-zero components in all $3$ factors. Since these
$3$ factors are not isomorphic unitary modules, we may conclude that
$\mathcal{L}$ is in fact a surjective map from $\ker(K_J)$ to $W_1\oplus W_2\oplus W_3$ as desired.\hfill\qedbox

\section{Para-Kaehler geometry}\label{sect-3}
The proof that every para-Kaehler curvature model is geometrically realizable by a para-Kaehler manifold is
essentially the same as the proof in the Kaehler setting given above in Section \ref{sect-2}.
Let $\tilde J$ be the  canonical integrable almost para-complex structure on $\mathbb{R}^m$ given by
$$\tilde J\partial_{1+2k}=\partial_{2+2k}\quad\text{and}\quad
\tilde J\partial_{2+2k}=\partial_{1+2k}\text{ for }0\le k<\bar m\,.$$
Let $g_0$ be the canonical flat para-Hermitian metric on $\mathbb{R}^m$;
$$ \tilde g_{0,ij}=\left\{\begin{array}{rl}
1&\text{if }i=j\equiv0\mod 2,\\
-1&\text{if }i=j\equiv1\mod2,\\
0&\text{otherwise}\,.\end{array}\right.$$
Again, we identify $(V,\langle\cdot,\cdot\rangle,\tilde J)$ with $(\mathbb{R}^m, \tilde g_0,\tilde J)$. Let $S_-^2(V^*)$ be the
space of symmetric
$2$-tensors
$\theta$ so
$\tilde J^*\theta=-\theta$. Given
$\Theta$ in
$S^2_-(V^*)\otimes S^2(V^*)$, we construct the para-Hermitian metric:
$${\tilde g_{\Theta,ij}}:= \tilde g_{0,ij}+ \Theta_{ijkl}u^ku^l\,.$$
We use Equation (\ref{eqn-2.c}) to define $K_{\tilde J}:S^2_-(V^*)\otimes S^2(V^*)\rightarrow\Lambda^3(V^*)\otimes V^*$ by
$$
\{(K_{\tilde J}\Theta)(x,y,z)\}(w):=\Theta(x,\tilde Jy,z,w)
+\Theta(y,\tilde Jz,x,w)+\Theta(z,\tilde Jx,y,w)\,.
$$
The curvature of  $\tilde g_\Theta$ at the origin is given by$$
({\mathcal{L}}\Theta)(x,y,z,w):=\Theta(x,z,y,w)+\Theta(y,w,x,z)-\Theta(x,w,y,z)-\Theta(y,z,x,w)\,.
$$

We may decompose $\mathfrak{A}$ as a direct sum of irreducible factors under the action of the para-unitary group. This decomposes the para-Kaehler tensors as a direct sum $\WWW_1\oplus\WWW_2\oplus\WWW_3$.
We show $\mathcal{L}$ is a surjective map from
$\ker(K_{\tilde J})$ to
$\WWW_1\oplus\WWW_2\oplus\WWW_3$ by taking
\begin{eqnarray*}
\Theta&=& \frac12(e^1\otimes e^1-e^2\otimes e^2)
\otimes(e^1\otimes e^1-e^2\otimes e^2),\\
\tilde g_\Theta&=& \tilde g_0
+ \frac12(u_1^2-u_2^2)(du_1^2-du_2^2)\,.\end{eqnarray*}
One now shows exactly as before that $A$ has non-zero components in each of these 3 factors; we omit details as the analysis is exactly the same as in the Kaehler setting. This shows every para-Kaehler curvature model is geometrically realizable by a
para-Kaehler manifold.\hfill\qed

\section{The Cauchy-Kovalevskaya Theorem}\label{sect-4}
In this section, we show the realization of Theorem \ref{thm-1.6} can be chosen to have constant scalar curvature; the
corresponding argument in the para-Kaehler case is similar and is therefore omitted. We shall adapt an argument given
in \cite{BGSV09} and refer to that paper for additional details. The major difference is that we consider a
$4^{\operatorname{th}}$ order quasi-linear partial differential equation rather than a $2^{\operatorname{nd}}$ order
equation.

We begin by
recalling the classical Cauchy-Kovalevskaya Theorem as formulated by Evans
\cite{E}. Set
$u=(y,u_m)$ for $y=(u_1,...,u_{m-1})\in\mathbb{R}^{m-1}$. Let $\Phi=\Phi(u)$ be real analytic. We consider the
$3^{\operatorname{rd}}$ order jets of $\Phi$:
$$\xi:=\{\Phi,\ \partial_{j_1}\Phi,\ \partial_{j_1}\partial_{j_2}\Phi,\
\partial_{j_1}\partial_{j_2}\partial_{j_3}\Phi\}\,.$$
We consider a $4^{\operatorname{th}}$ order quasi-linear equation in $\Phi$:
\begin{equation}\label{eqn-4.a}
\psi^{i_1i_2i_3i_4}(\xi)\partial_{i_1}\partial_{i_2}\partial_{i_3}\partial_{i_4}\Phi+
\psi(\xi)=0\end{equation}
where the coefficients $\psi^{i_1i_2i_3i_4}$ and $\psi$ are real analytic functions of the variables $\xi$.
Impose the Cauchy data $\xi=0$ on the initial hypersurface  $u_m=0$, i.e.
\begin{equation}\label{eqn-4.b}
\Phi(y,0)=0,\quad\partial_m\Phi(y,0)=0,\quad\partial_m\partial_m\Phi(0,y)=0,\quad
\partial_m\partial_m\partial_m\Phi(y,0)=0\,.
\end{equation}
\begin{theorem}\label{thm-4.1}
If $\psi^{mmmm}(0)\ne0$, there is $\epsilon>0$ and a unique real analytic $\Phi$ defined
for
$|u|<\epsilon$ which satisfies Equations {\rm(\ref{eqn-4.a})} and {\rm(\ref{eqn-4.b})}.
\end{theorem}

We use Theorem \ref{thm-4.1} to show the metric of Theorem \ref{thm-1.6} can be chosen to be real analytic and
with constant scalar curvature. We adopt the notation of Section \ref{sect-2}. For $0\le k<\bar m$, set:
$$\begin{array}{ll}
z_k:=u_{1+2k}+\sqrt{-1}u_{2+2k},&\bar z_k:=u_{1+2k}-\sqrt{-1}u_{2+2k},\\
dz_k:=du_{1+2k}+\sqrt{-1}du_{2+2k},&d\bar z_k:=du_{1+2k}-\sqrt{-1}du_{2+2k},\vphantom{\vrule height 11pt}\\
\partial_{z_k}:=\textstyle\frac12\{\partial_{1+2k}-\sqrt{-1}\partial_{2+2k}\},&
\partial_{\bar z_k}:=\textstyle\frac12\{\partial_{1+2k}+\sqrt{-1}\partial_{2+2k}\}\,.\vphantom{\vrule height 11pt}
\end{array}$$

 Let $\circ$ denote the symmetric product. If $\Phi$ is a real analytic function (which is called the {\it Kaehler potential}), form
$$\kappa_\Phi:=\left\{\partial_{z_j}\partial_{\bar z_k}\Phi\right\}dz^j\circ d\bar z^k\in S^2_+(V^*)\,.$$
This is a $J$-invariant symmetric real $2$-tensor with $d\Omega_{\kappa_\Phi}=0$.

Suppose given a Kaehler curvature model $\mathfrak{C}=(V,\langle\cdot,\cdot\rangle,J,A)$
with scalar curvature
$c$. Use the methods of Section \ref{sect-2} to choose
$\Theta$ so that the curvature tensor of the metric $g_\Theta$ is given by $A$ at the origin where
$g_\Theta$ is the real analytic metric defined by Equation (\ref{eqn-2.b}) and where $J$ is the integrable
complex structure given by Equation (\ref{eqn-2.a}). Consider the metric
$h_{\Theta,\Phi}:=g_\Theta+\kappa_\Phi$. We impose the Cauchy initial data given by Equation (\ref{eqn-4.b}) so
$h_{\Theta,\Phi}(0)=\delta$ and thus $h_{\Theta,\Phi}$ is a Riemannian metric in some neighborhood of the origin. The
scalar curvature
$\tau_\Phi$ of $h_{\Theta,\Phi}$ is given by a quasi-linear $4^{\operatorname{th}}$-order equation of the form given by
Equation (\ref{eqn-4.a}).

We suppress terms which do not involve maximal derivatives in $\partial_m$ to write:
$$
 h_{\bar m\bar m}=h_{mm}=\textstyle\frac14\partial_m^2\Phi+...,\quad
R_{m\bar m\bar m m}=-\textstyle\frac18\partial_m^4\Phi+...,\quad
 \tau_\Phi=-\textstyle\frac14\partial_m^4\Phi+...\,.
$$
Thus the non-degeneracy condition of Theorem \ref{thm-4.1} is
satisfied and we can solve the equation $\tau_\Phi-c=0$ with vanishing Cauchy initial data. All the
$4^{\operatorname{th}}$ order derivatives of
$\Phi$ vanish except possibly for
$\partial_m^4\Phi$ -- the equation $\tau_\Phi-c=0$ implies this vanishes as well. Thus $\Phi=O(|u|^5)$ so $\Phi$
makes no contribution to the curvature tensor at the origin. This shows the Kaehler manifold in question
can be chosen to have constant scalar curvature; a similar calculation using a para-Kaehler potential
pertains in the para-Kaehler setting as well. This completes the proof of all the assertions of this
paper.\hfill\qedbox

\section*{Acknowledgments}

Research of M. Brozos-V\'azquez partially supported by Project MTM2006-01432 (Spain).
Research of P. Gilkey partially supported by DFG PI 158/4-6 (Germany) and by Project MTM2006-01432 (Spain).
Research of E. Merino partially supported by FEDER and Project MTM2008-05861 MICINN (Spain).


\begin{thebibliography}{AAA}

\bibitem{AD01} V. Apostolov, and T. Dra\v ghici, {\it The curvature and the integrability of almost-K\"ahler manifolds: a survey}.  Symplectic and contact topology: interactions and perspectives (Toronto, ON/Montreal, QC, 2001),  25--53, Fields Inst. Commun., 35, Amer. Math. Soc., Providence, RI, 2003.

\bibitem{B08} I. Biswas,
{\it Holomorphic principal bundles with an elliptic curve as the structure group},
Int. J. Geom. Methods Mod. Phys. {\bf 5} (2008), 851--862.

\bibitem{BGNS08} N. Blazic, P. Gilkey, S. Nik\v cevi\'c, and I. Stavrov,
{\it Curvature structure of self-dual 4-manifolds}, Int. J. Geom. Methods Mod. Phys. {\bf 5} (2008), 1191-1204.

\bibitem{B07}
R. Bos, {\it Geometric quantization of Hamiltonian actions of Lie algebroids and Lie groupoids},  Int. J. Geom. Methods Mod. Phys.  {\bf 4}
(2007), 389--436.

\bibitem{BGKSW09} M. Brozos-V\'azquez, P. Gilkey, H. Kang, S. Nik\v cevi\'c, and G. Weingart,
{\it Geometric realizations of curvature models by manifolds with constant scalar curvature},
Differential Geom. Appl., doi:10.1016/j.difgeo.2009.05.002  (arXiv:0811.1651).

\bibitem{BGKS09} M. Brozos-V\'azquez, P. Gilkey, H. Kang, and S. Nik\v cevi\'c,
{\it Geometric realizations of Hermitian curvature models}, J. Math. Soc. Japan, to appear (arXiv:0812.2743).

\bibitem{BGSV09} M. Brozos-V\'azquez, P. Gilkey, S. Nik\v cevi\'c, and R. V\'{a}zquez-Lorenzo, {\it Geometric Realizations of
para Hermitian curvature models}, Results in Math., to appear (arXiv:0902.1697).

\bibitem{CM08} J. Clemente-Gallardo, and G. Marmo, {\it Basics of quantum mechanics, geometrization and some applications to quantum information},  Int. J. Geom. Methods Mod. Phys.  {\bf 5}  (2008),  989--1032.

\bibitem{CVLS06} V. Cort\'es, M. Lawn, and L. Schaefer,
{\it Affine hyperspheres associated to special para-Kaehler manifolds}, Int. J. Geom. Methods Mod. Phys. {\bf 3} (2006), 995--1009.

\bibitem{CDG08} A. Cort\'es-Ayaso,  J. D\'i{a}z-Ramos, and E. Garc\'ia-R\'io, {\it Four-dimensional manifolds with degenerate self-dual Weyl curvature operator},  Ann. Global Anal. Geom.  {\bf 34} (2008), 185--193.

\bibitem{DJ08} L. Dallagnol, and M. Jardim, {\it Nonsingular complex instantons on Euclidean spacetime},   Int. J. Geom. Methods Mod. Phys.  {\bf 5}  (2008),  963--971.

\bibitem{DDGMMV07} J. Davidov, J. D\'i{a}z-Ramos, E. Garc\'ia-R\'io, Y. Matsushita, O. Muskarov, and V\'azquez-Lorenzo, {\it    Almost Kaehler Walker 4-manifolds},  J. Geom. Phys.  {\bf57}  (2007),  1075--1088.

\bibitem{E} L. Evans, {\bf Partial Differential Equations}, Graduate Texts in Mathematics {\bf 19}, American
Mathematical Society, Providence R. I.

\bibitem{GHVV08}    E. Garc\'ia-R\'io, A. Haji-Badali, M. V\'azquez-Abal, and R. V\'azquez-Lorenzo, {\it Lorentzian
3-manifolds with commuting curvature operators},  Int. J. Geom. Methods Mod. Phys.  {\bf5} (2008), 557--572.

\bibitem{GN07} P. Gilkey, and S. Nik\v cevi\'c, {\it Pseudo-Riemannian Jacobi-Videv manifolds},  Int. J. Geom. Methods Mod. Phys.  {\bf4}  (2007),  727--738.
\bibitem{GN08x} P. Gilkey, and S. Nik\v cevi\'c, {\it     Geometrical representations of equiaffine curvature operators},  Results in Mathematics {\bf52} (2008), 281--287.

 \bibitem{GNW09x} P. Gilkey, S. Nik\v cevi\'c, and D. Westerman,  {\it   Geometric realizations of generalized algebraic curvature operators}, J. Math. Phys.
{\bf50} (2009), 013515.

 \bibitem{GNW09} P. Gilkey, S. Nik\v cevi\'c, and D. Westerman, {\it Riemannian geometric realizations for Ricci
tensors of generalized algebraic curvature operators}, Differential Geometry - Proceedings of
the VIII International Colloquium, Santiago de Compostela, Spain. Ed. J. Alvarez-Lopez and E. Garcia-Rio. World Scientific (2009),
175--184.

\bibitem{gray} A. Gray, {\it Curvature identities for Hermitian and
almost Hermitian manifolds}, {T{\^o}hoku Math. J.} {\bf 28}
(1976), 601--612.

\bibitem{K08} H. Kang,
{\it Compactness of $D$-isospectral metrics},
Int. J. Geom. Methods Mod. Phys. {\bf5} (2008), 49--61.

\bibitem{MP09} R. Milson, and N. Pelavas, {\it The curvature homogeneity bound for Lorentzian four-manifolds}, Int. J.
Geom. Methods Mod. Phys. {\bf 6} (2009), 99--127.

\bibitem{NN57}
A. Newlander, and L. Nirenberg,
Complex analytic coordinates in almost complex manifolds,
{\it Ann. of Math.} {\bf 65} (1957), 391--404.


\bibitem{T08} R. Tresguerres,  {\it Dynamically broken anti-de Sitter action for gravity},  Int. J. Geom. Methods Mod. Phys. {\bf 5}  (2008), 171--183.

\bibitem{TV81} F. Tricerri, and L. Vanhecke, {\it Curvature tensors on almost Hermitian manifolds}, {Trans.
Amer. Math. Soc.} {\bf 267} (1981), 365--397.

\bibitem{W08} S. Watterson, {\it The chiral and flavor projection of Dirac-Kaehler fermions in the geometric discretization},  Int. J. Geom. Methods Mod. Phys.  {\bf 5}  (2008),  345--362.
\end{thebibliography}
\end{document}